\documentclass[12pt]{article}
\usepackage{amssymb}
\usepackage{amsmath}
\usepackage{listings}
\usepackage{color}
\lstdefinestyle{customc}{
  belowcaptionskip=1\baselineskip,
  breaklines=true,
  frame=single,
  xleftmargin=\parindent,
  language=C++,
  showstringspaces=false,
  basicstyle=\footnotesize\ttfamily,
  keywordstyle=\bfseries\color{red},
  commentstyle=\itshape\color{magenta},
  identifierstyle=\color{blue},
}
\newtheorem{teo}{Theorem}
\newtheorem{prop}{Proposition}

\newtheorem{cor}{Corollary}
\newtheorem{rem}{Remark}

\newenvironment{Proof} 
{\par\noindent{\bf Proof.}} 
{\hfill$\scriptstyle\blacksquare$}

\title{ On refinement masks of tight wavelet frames
}
\author{E. A. Lebedeva\footnote{Mathematics and Mechanics Faculty, Saint Petersburg State University,
Universitetsky prospekt, 28, Peterhof,  Saint Petersburg, 198504, Russia}, 
 I. A. Shcherbakov\footnote{Mathematics and Mechanics Faculty, Saint Petersburg State University,
Universitetsky prospekt, 28, Peterhof,  Saint Petersburg, 198504, Russia}
}
\date{
 ealebedeva2004@gmail.com, stscherbakov99@yandex.ru
}

\begin{document}
\maketitle

\newcommand{\nul}{{\bf0}}
\newcommand{\rd}{{\mathbb R}^d}
\newcommand{\zd}{{\mathbb Z}^{d}}
\renewcommand{\r}{{\mathbb R}}
\newcommand{\z} {{\mathbb Z}}
\newcommand{\cn} {{\mathbb C}}
\newcommand{\n} {{\mathbb N}}

\begin{abstract}
In the paper we obtain  sufficient conditions for a trigonometric polynomial to be a refinement mask corresponding to a tight wavelet frame. The condition is formulated in terms of the roots of a mask. In particular, it is proved that any trigonometric polynomial can serve as a mask if its associated algebraic polynomial has only negative roots (at least one of them, of course,  equals $-1$).               
\end{abstract}



\section{Introduction}

The unitary extension principle (UEP) of Ron and Shen \cite{RonShen} is one of the main tools for the constructions of tight wavelet frames. We recall it here for completeness of the presentation. 

\textbf{The unitary extension principle.} 
Let $\varphi \in L_2(\mathbb{R})$ be a refinement function, that is the following  refinement equation holds 
\begin{equation}
	\label{ref_eq}
	\widehat{\varphi}(\xi) = m_0(\xi/2) \widehat{\varphi}(\xi/2) \mbox{ a.e.,} 
\end{equation}
where $m_0 \in L_2(0,\,1)$ is a refinement mask. Let $\widehat{\varphi}$ be continuous at zero.  Let $m_1, \dots, m_r$ be $1$-periodic functions in $L_2(0,\,1)$ called wavelet masks such that the matrix
$$
M(\xi) = \left(
\begin{array}{llll}
m_0(\xi) & m_1(\xi) & \dots & m_r(\xi) \\
m_0(\xi+1/2) & m_1(\xi+1/2) & \dots & m_r(\xi+1/2) \\ 
\end{array}
\right)
$$                
satisfies the equality 
 \begin{equation}
	\label{unit_m}
	M(\xi) M^{\ast}(\xi) = I_2, \mbox{ a.e.,} 
\end{equation}
where $I_2$ is the identity matrix of size $2.$ Define in the Fourier domain wavelets $\psi^{(1)}, \dots, \psi^{(r)}$ by the formulas 
\begin{equation}
	\label{wav_eq}
	\widehat{\psi^{(k)}}(\xi) = m_k(\xi/2) \widehat{\varphi}(\xi/2) \mbox{ a.e., } k=1,\dots,r  
\end{equation}
 Then  the  system  of  functions $\psi_{j,k}^{(k)}$, $j,k \in \mathbb{Z},$ $k=1,\dots,r,$   forms  a  tight  frame  in  $L_2(\mathbb{R})$  with frame bounds 
$\displaystyle A=B=\left|\widehat{\varphi}(0)\right|^2.$

It is well known \cite{Petukhov} that the general setup together with the inequality 
\begin{equation}
	\label{mask_in}
	\left|m_0(\xi)\right|^2 + \left|m_0(\xi+1/2)\right|^2 \leq 1 \mbox{ a.e.} 
\end{equation}
always provides a solution for matrix equation (\ref{unit_m}) and  makes it possible to obtain a  frame with two wavelet generators $\psi^{(1)},$ $\psi^{(2)}$. So, if we wish to design a frame by means of UEP we need to find a function $\varphi \in L_2(\mathbb{R})$ such that $\widehat{\varphi}$ is continuous at zero, it satisfies (\ref{ref_eq}) and its mask satisfies (\ref{mask_in}). 
One can construct a wavelet frame starting with a refinement mask. In this case a refinement function is defined by an infinite product $\displaystyle \widehat{\varphi}(\xi) = \prod_{j=1}^{\infty} m_0(\xi/2^j)$ and one needs to check that $\varphi \in L_2(\mathbb{R}).$ It can be done using  

\textbf{the Mallat theorem} \cite[Lemma 4.1.3]{NPS}. Suppose $m_0(\xi) = \sum_{k\in\mathbb{Z}} c_k {\rm e}^{2\pi {\rm i} k \xi},$ $m_0(0)=1,$ $c_k = O(|k|^{-2-\varepsilon}),$ $\varepsilon>0,$ and (\ref{mask_in}) holds. Define $\displaystyle \widehat{\varphi}(\xi) = \prod_{j=1}^{\infty} m_0(\xi/2^j)$. Then 
$\varphi \in L_2(\mathbb{R}),$ and $\|\varphi\|\leq 1.$

Let $m_0$ be a trigonometric polynomial and $m_0(0)=1$. Then it is well known that  the function  $\displaystyle \widehat{\varphi}(\xi) = \prod_{j=1}^{\infty} m_0(\xi/2^j)$ is an entire function of exponential  type, thus $\widehat{\varphi}$ is continuous at zero. 
If additionally $m_0$ satisfies inequality  (\ref{mask_in}),  then according to the Mallat theorem the function  $\displaystyle \widehat{\varphi}(\xi) = \prod_{j=1}^{\infty} m_0(\xi/2^j) \in L_2(\mathbb{R})$, it is a corresponding refinement function, and it generates a tight wavelet frame as it is provided by UEP. 

The purpose of the paper is to obtain  sufficient conditions for roots of a trigonometric polynomial 
$m_0$ to   satisfy   (\ref{mask_in}). As it is seen from the above, inequality (\ref{mask_in}) is a cornerstone  for the constructions of tight wavelet frames. It is worth to note that in the case of orthogonal wavelets assumptions on a refinement mask are much more restrictive 
(see \cite[Theorem 4.1.2]{NPS}). The paper is organized as follows. In Section 2, we first consider polynomials of low degrees ($2$ and $3$) and obtain   not only a sufficient, but also a necessary condition to satisfy (\ref{mask_in}). This is done in Proposition 1 and Proposition 2. Then in Theorem \ref{main_th}, we consider polynomials of a arbitrary degree and obtain a sufficient condition to satisfy (\ref{mask_in}).
A particular, extremely easily checked case of Theorem \ref{main_th} is formulated in Corollary 2. In Section 3, we present a validation algorithm for the sufficient condition obtained in Theorem   \ref{main_th}.  

\section{Results}
\subsection{Preliminary}
Consider the algebraic polynomial $P(z)$ associated with $m_0(\xi)$, which is given by the equation $m_0(\xi)=P(e^{2\pi i\xi})e^{2\pi i\xi \max\{-N,0\}}.$ It is easy to see that the inequality

\begin{equation}
	\label{mask_in2}
	\left|P(z)\right|^2 + \left|P(-z)\right|^2 \leq 1 \mbox{ a.e. on $\mathbb{T}$}
\end{equation}
is equivalent (\ref{mask_in}) and $m_0(0)=1$ iff $P(1)=1.$

We immediately make a couple of obvious remarks.

\begin{rem}
If an algebraic polynomial $P(z)$ satisfies the inequality (\ref{mask_in2}) and $P(1)=1,$ then this polynomial has a root at the point $z=-1$.
\end{rem}

\begin{rem}\label{r2}
If $P(1)=1$, then the polynomial $P$ can be written in the form $\displaystyle P(z) = \prod\limits_{i=1}^n\frac{z-z_i}{1-z_i},$ where $z_i\not = 1.$
\end{rem}

\begin{rem}
The polynomial $P(z)$ of degree at least $2$ satisfies (\ref{mask_in2}) and $P(1)=1$ iff the inequality (\ref{mask_in2}) holds for $Q(z)=zP(z)$.
\end{rem}

Let the function $\psi_{z_0}$ be given by $z\in\mathbb{T} \rightarrow \psi_{z_0}(z) = \Bigl|\dfrac{z-z_0}{1-z_0}\Bigr|^2,$ where $z_0\not=1$. If polynomial $P(z)$ has roots $z_1, \dots, z_n$ and $P(1)=1,$ then according to Remark \ref{r2} $|P(z)|^2=\prod\limits_{i=1}^n\psi_{z_i}(z)$. Suppose $\alpha\in Arg(z-z_0)$ and $\beta\in Arg(z_0)$, then by the cosine theorem it follows that 
$$\psi_{z_0}(z)=\frac{1+|z_0|^2-2|z_0|\cos \alpha}{1+|z_0|^2-2|z_0|\cos \beta},$$
$$\psi_{z_0}(e^{i\varphi})=\frac{1+|z_0|^2-2|z_0|\cos (\varphi-\beta)}{1+|z_0|^2-2|z_0|\cos \beta}$$

Denote by $x={\rm Re} \, z_0$ and $y={ \rm  Im }\, z_0$, we get a trigonometric first-degree polynomial depending on  $x,y$ 
$$\psi_{z_0}(e^{i\varphi})=\frac{1+x^2+y^2-2x\cos \varphi-2y\sin \varphi}{1+x^2+y^2-2x}=$$
$$= 1 + \frac{2x}{(x-1)^2+y^2}-\frac{2x}{(x-1)^2+y^2}\cos \varphi-\frac{2y}{(x-1)^2+y^2}\sin \varphi.$$

Denote  by 
 \begin{gather}
  \notag
     F_1(x,y)=1 + \frac{2x}{(x-1)^2+y^2},\\
    \label{f1-3}
 F_2(x,y)=-\frac{2x}{(x-1)^2+y^2}, \\ 
 \notag
 F_3(x,y)=-\frac{2y}{(x-1)^2+y^2}.
 \end{gather}
  the coefficients of this polynomial, then

$$\psi_{z_0}(z)=\psi_{z_0}(e^{i\varphi})=F_1(x,y)+F_2(x,y)\cos \varphi+F_3(x,y)\sin \varphi,$$
and 
$$\psi_{z_0}(-z)=\psi_{z_0}(e^{i(\pi+\varphi)})=F_1(x,y)-F_2(x,y)\cos \varphi-F_3(x,y)\sin \varphi.$$
So, the left-hand side of inequality (\ref{mask_in}) takes the form
\begin{equation}
\label{def_T}
    |P(z)|^2+|P(-z)|^2=\prod_{i=1}^n\psi_{z_i}(z)+\prod_{i=1}^n\psi_{z_i}(-z)
    =:T(\varphi).
\end{equation}

In the sequel, to obtain sufficient conditions for the roots $z_1, \dots, z_n$ of the polynomial $P(z)$ to satisfy inequality (\ref{mask_in}), we study the trigonometric polynomial $T(\varphi)$. Note that $T$  has degree at most $n$, it is nonnegative, $\pi$-periodic. So, it only contains monomials with even angles.

\subsection{The simplest cases}
First, assume that the degree of the polynomial $P$ is equal to $2$. In this case we obtain not only a sufficient, but also a necessary condition to satisfy (\ref{mask_in}). 

\begin{prop}
Let $P(z)$ be an algebraic polynomial and it satisfies a condition $P(1)=1$. The numbers $z_1 = -1$ and $z_2$ are all roots of this polynomial. Then (\ref{mask_in2}) holds if and only if  $z_2\le 0.$ 
\end{prop}
\begin{Proof} 
 In this case $F_1(z_1)=F_2(z_1)=\dfrac 1{2}$ and $F_3(z_1)=0$. Denote by $A=F_1(z_2)$ and $B=F_3(z_2)$. Substituting into (\ref{def_T}),  we get 

$$T(\varphi)=\Bigl(\dfrac 1{2}+\dfrac 1{2}\cos\varphi\Bigr)(A+(1-A)\cos\varphi+B\sin\varphi)+\Bigl(\dfrac 1{2}-\dfrac 1{2}\cos\varphi\Bigr)(A-(1-A)\cos\varphi-B\sin\varphi)=$$
$$=A+(1-A)\cos^2\varphi+B\sin\varphi\cos\varphi=\dfrac {1+A}2+\dfrac {1-A}2\cos 2\varphi +\dfrac B{2}\sin 2\varphi\le
$$
$$
\le \dfrac {1+A}2+\dfrac {1}2\sqrt{(1-A)^2+B^2}.$$

Since this inequality is exact, it follows that inequality (\ref{def_T}) is equivalent to 
$$\dfrac {1+A}2+\dfrac {1}2\sqrt{(1-A)^2+B^2} \le 1,$$
that is
$$\sqrt{(1-A)^2+B^2} \le 1-A.$$
That is equivalent to two conditions
$$
\begin{cases}
1-A\ge 0 \Leftrightarrow F_2(z_2)\ge 0 \Leftrightarrow {\rm Re} \ z_2 \le  0, \\
B=0 \Leftrightarrow F_3(z_2)= 0 \Leftrightarrow {\rm Im} \ z_2 =  0. 
\end{cases}
$$
\end{Proof}

Similarly we get the necessary and sufficient conditions for $n = 3.$

\begin{prop} \label{degree3}
Let $P(z)$ be an algebraic polynomial and it satisfies the condition $P(1)=1$. The numbers $\frac1{2}, \ z_1$ and $z_2$ are all roots of this polynomial. Denote $A_1=F_1(z_1), \ A_2=F_1(z_2), \ B_1=F_3(z_1)$ and $B_2=F_3(z_2).$ Then (\ref{mask_in2}) holds if and only if
$$\begin{cases}
1-A_1A_2-B_1B_2\ge 0, \\
B_1+B_2=0. 
\end{cases}$$
\end{prop}
\begin{Proof}
As in the previous proof, we substitute all the notations into (\ref{def_T}) and get 
\begin{multline*}
T(\varphi)=\Bigl(\dfrac 1{2}+\dfrac 1{2}\cos\varphi\Bigr)\Bigl(A_1+(1-A_1)\cos\varphi+B_1\sin\varphi\Bigr)\Bigl(A_2+(1-A_2)\cos\varphi+B_2\sin\varphi\Bigr)+ \\
+\Bigl(\dfrac 1{2}-\dfrac 1{2}\cos\varphi\Bigl)\Bigr(A_1-(1-A_1)\cos\varphi-B_1\sin\varphi\Bigr)\Bigl(A_2-(1-A_2)\cos\varphi-B_2\sin\varphi\Bigr)=
\end{multline*}
\begin{multline*}
=A_1A_2+\cos^2\varphi((1-A_1)(1-A_2)+(1-A_1)A_2+(1-A_2)A_1)+\sin^2\varphi\cdot B_1B_2+ \\
+\cos\varphi\sin\varphi((1-A_1)B_2+(1-A_2)B_1+A_1B_2+A_2B_1)=
\end{multline*}

$$=\dfrac{1+A_1A_2+B_1B_2}2+\cos 2\varphi \dfrac{1-A_1A_2-B_1B_2}2+\sin 2\varphi \dfrac{B_1+B_2}2\le$$

$$\le \dfrac{1+A_1A_2+B_1B_2}2+\sqrt{\dfrac{(1-A_1A_2-B_1B_2)^2}4+\dfrac{(B_1+B_2)^2}4}.$$
Since this inequality is exact, it follows that inequality (\ref{def_T}) is equivalent to 
$$\dfrac{1+A_1A_2+B_1B_2}2+\sqrt{\dfrac{(1-A_1A_2-B_1B_2)^2}4+\dfrac{(B_1+B_2)^2}4}\le 1.$$
It can be rewritten as
$$\sqrt{(1-A_1A_2-B_1B_2)^2+(B_1+B_2)^2}\le1-A_1A_2-B_1B_2,$$
that is equivalent to
$$
\begin{cases}
1-A_1A_2-B_1B_2\ge 0, \\
B_1+B_2=0.
\end{cases}
$$
\end{Proof}

\begin{cor}
Let $P(z)$ be an algebraic polynomial of degree 3 with real roots $x_1, x_2$ and $-1$. It satisfies the condition $P(1)=1$. Then (\ref{mask_in2}) holds if and only if $$x_1x_2(x_1+x_2-2)+x_1+x_2\le0$$
\end{cor}
\begin{Proof}
Note that for real roots we get $B_1=B_2=0,$ so the necessary and sufficient conditions from  Proposition \ref{degree3} can be written as follows

$$\Bigl(1+\dfrac{2x_1}{(x_1-1)^2}\Bigr)\Bigl(1+\dfrac{2x_2}{(x_2-1)^2}\Bigr)\le 1
\Leftrightarrow
4x_1x_2+2x_1(x_2-1)^2+2x_2(x_1-1)^2\le 0$$
$$
\Leftrightarrow
x_1x_2(x_1+x_2-2)+x_1+x_2\le0$$
\end{Proof}

\subsection{The main result}
Now we return to polynomials of degree $n$. We consider a trigonometrical polynomial $m_0$ and  the algebraic polynomial $P$ associated with $m_0$, that is $m_0(\xi)=P(e^{2\pi i\xi})e^{2\pi i\xi \max\{-N,0\}}.$
Suppose all  roots $z_1, \dots, z_n$ of the polynomial $P$ are real.   Denote by $a_i=F_1(z_i,0),$ $i\in[1..n]$ (see (\ref{f1-3})). Denote by $\sigma_k$ the elementary symmetric polynomials $\displaystyle \sum\limits_{\substack{S\subset [1..n]\\ \#S=k}}\prod_{j\in S}a_j$ for all $k\in[0..n]$. Put $\displaystyle \rho_k:=\dfrac{\sigma_k}{{n \choose k}}$  as \cite[page 73]{Polynom}. Now we are ready to formulate the main theorem.

\begin{teo}
\label{main_th}
Let $T(\varphi)$ be a trigonometric polynomial  built by the polynomial $P$ as it is defined in (\ref{def_T}).  Let all roots $z_1, \dots, z_n$ of the polynomial $P$ be real 
 and at least one of them equals $-1$. Then if 
$$\Delta^{2k}\rho_{n-2k}=\sum_{j=n-2k}^n{2k \choose n-j}(-1)^{n-j}\rho_{j}\ge0$$
for any $k\in[0..[\frac n{2}]],$ then $T(\varphi)\le 1$ for any $\varphi\in\mathbb{R}.$
\end{teo}
\begin{Proof}
Taking into account that for $z_i \in \r$, $i\in[1..n]$, we get $F_2(z_i,0) =1- a_i$, and $F_3(z_i,0) = 0$  (see (\ref{f1-3})), the function $\psi_{z_i}$ has a simpler form

$$\prod_{i=1}^n\psi_{z_i}(e^{i\varphi})=\prod_{i=1}^n(a_i+(1-a_i)\cos \varphi)=\sum_{k=0}^n\cos^k \varphi\cdot\sum\limits_{\substack{S\subset [1..n]\\ \#S=k}}\prod_{j\in S}(1-a_j)\prod_{j\not\in S}a_j.$$
Therefore, we have
\begin{equation}
\label{T1}
    T(\varphi)=\prod_{i=1}^n\psi_{z_i}(e^{i\varphi})+\prod_{i=1}^n\psi_{z_i}(e^{i(\pi+\varphi)})=2\sum_{k=0}^{[\frac n{2}]}\cos^{2k} \varphi\cdot\sum\limits_{\substack{S\subset [1..n]\\ \#S=2k}}\prod_{j\in S}(1-a_j)\prod_{j\not\in S}a_j.
\end{equation}

Firstly,  we use Euler's formula for a cosine with $k\ge1$

$$\cos^{2k} \varphi=\Bigl(\frac{e^{i\varphi}+e^{-i\varphi}}{2}\Bigr)^{2k}=\frac 1{2^{2k}}\Bigl(\sum_{l=0}^{2k}{2k \choose l}e^{i\varphi l}e^{-i\varphi (2k-l)}\Bigr)=$$
$$=\frac 1{2^{2k}}\Bigl({2k \choose k}+\sum_{l=0}^{k-1}{2k \choose l}(e^{2i\varphi (l-k)}+e^{2i\varphi (k-l)})\Bigr)=
\frac 1{2^{2k}}\Bigl({2k \choose k}+2\sum_{l=0}^{k-1}{2k \choose l}\cos2(l-k)\varphi\Bigr)=$$
$$=\frac 1{2^{2k}}\Bigl({2k \choose k}+2\sum_{l=1}^{k}{2k \choose k-l}\cos2l\varphi\Bigr).$$

Secondly, we consider the coefficients of $T(\varphi)$ as combinations of symmetric polynomials of $a_1, \dots, a_n$:

$$\sum\limits_{\substack{S\subset [1..n]\\ \#S=2k}}\prod_{j\in S}(1-a_j)\prod_{j\not\in S}a_j=
\sum\limits_{\substack{S\subset [1..n]\\ \#S=2k}}\Bigr(\prod_{j\not\in S}a_j\cdot\sum_{l=0}^{2k}\Bigr((-1)^l\sum\limits_{\substack{T\subset S\\ \#T=l}}\prod_{m\in T}a_m\Bigl)\Bigl)=$$
$$=\sum\limits_{\substack{S\subset [1..n]\\ \#S=2k}}\sum_{l=0}^{2k}\Bigr((-1)^l\sum\limits_{\substack{T\subset S\\ \#T=l}}\prod_{m\in T\cup S^C}a_m\Bigl)=\sum_{l=0}^{2k}\Bigr((-1)^l\sum\limits_{\substack{S\subset [1..n]\\ \#S=2k}}\sum\limits_{\substack{T\subset S\\ \#T=l}}\prod_{m\in T\cup S^C}a_m\Bigl)=$$
$$=\sum_{l=0}^{2k}\Bigr((-1)^l\sum\limits_{\substack{S\subset [1..n]\\ \#S=n-2k}}\sum\limits_{\substack{T\cap S = \emptyset\\ \#T=l}}\prod_{m\in T\cup S}a_m\Bigl)=
\sum_{l=0}^{2k}\Bigr((-1)^l{n-2k+l \choose l}\sum\limits_{\substack{S\subset [1..n]\\ \#S=n-2k+l}}\prod_{m\in S}a_m\Bigl)=$$
$$=\sum_{l=0}^{2k}\Bigr((-1)^l{n-2k+l \choose l}\sigma_{n-2k+l}\Bigl)=\sum_{l=0}^{2k}\Bigr((-1)^l{n-l \choose 2k-l}\sigma_{n-l}\Bigl)$$

Thus, 

$$T(\varphi)=2\sigma_{n}+2\sum_{k=1}^{[\frac n{2}]}\Bigl(\frac1{2^{2k}}\Bigl({2k \choose k}+2\sum_{l=1}^{k}{2k \choose k-l}\cos2l\varphi\Bigr)\cdot\sum_{m=0}^{2k}\Bigl((-1)^m{n-m \choose 2k-m}\sigma_{n-m}\Bigr)\Bigr)=$$

\begin{multline*}
=2\sigma_{n}+2\sum_{k=1}^{[\frac n{2}]}\sum_{m=0}^{2k}\frac 1{2^{2k}}{2k \choose k}(-1)^m{n-m \choose 2k-m}\sigma_{n-m}
+   \\
+4\sum_{k=1}^{[\frac n{2}]}\sum_{l=1}^{k}\sum_{m=0}^{2k}(-1)^m\frac 1{2^{2k}}{2k \choose k-l}\cos2l\varphi {n-m \choose 2k-m}\sigma_{n-m}=
\end{multline*}

\begin{multline*}
=2\sum_{k=0}^{[\frac n{2}]}\sum_{m=0}^{2k}\frac 1{2^{2k}}{2k \choose k}(-1)^m{n-m \choose 2k-m}\sigma_{n-m}
+   \\
+4\sum_{l=1}^{[\frac n{2}]}\cos2l\varphi\Bigl(\sum_{k=l}^{[\frac n{2}]}\sum_{m=0}^{2k}(-1)^m\frac 1{2^{2k}}{2k \choose k-l} {n-m \choose 2k-m}\sigma_{n-m}\Bigr)
\end{multline*}

To simplify notations, we denote the  coefficients by 
$$ 
d_l=4\sum\limits_{k=l}^{[\frac n{2}]}\sum\limits_{m=0}^{2k}(-1)^m\frac 1{2^{2k}}{2k \choose k-l} {n-m \choose 2k-m}\sigma_{n-m},
$$
where $l\in[0,\dots,[\frac n{2}]]$.
Then $T(\varphi)$ takes the form $\displaystyle  T(\varphi)=\frac {d_0}{2}+\sum\limits_{l=1}^{[\frac n{2}]}d_l\cos2l\varphi.$ Note that

$$d_l=4\sum\limits_{k=l}^{[\frac n{2}]}\frac 1{2^{2k}}\sum\limits_{m=0}^{2k}(-1)^m{2k \choose k-l} {n-m \choose 2k-m}{n \choose n-m}\rho_{n-m}=$$
$$=4\sum\limits_{k=l}^{[\frac n{2}]}\frac 1{2^{2k}}{n \choose n-2k,k-l,k+l}\sum\limits_{m=0}^{2k}(-1)^m{2k \choose m}\rho_{n-m}=$$
$$=4\sum\limits_{k=l}^{[\frac n{2}]}\frac 1{2^{2k}}{n \choose n-2k,k-l,k+l}\Delta^{2k}\rho_{n-2k}\ge0$$

Therefore, all coefficients of $T(\varphi)$ are nonnegative and $T(\varphi)\le T(0)=\prod\limits_{i=1}^n\psi_{z_i}(1)+\prod\limits_{i=1}^n\psi_{z_i}(-1)=1+\prod\limits_{i=1}^n(2a_i-1).$ Since at least one of $z_1, \dots, z_n$ equals $-1$, it follows that at least one of $a_1, \dots, a_n$ equals $F_1(-1)=1/2$. Therefore, $T(0)=1$. As was to be proved.
\end{Proof}

\begin{cor}
Let $T(\varphi)$ be a trigonometric polynomial  built by the polynomial $P$ as it is defined in (\ref{def_T}).  Let all  roots $z_1, \dots, z_n$ of the polynomial $P$ be less or equal to zero  
 and at least one of them equals $-1$.
 Then $T(\varphi)\le 1$ for any $\varphi\in\mathbb{R}$.
\end{cor}

\begin{Proof}
It is clear, that $z_i\le 0 \Leftrightarrow a_i\le 1.$ Then using (\ref{T1}) we get
$$T(\varphi)=2\sum_{k=0}^{[\frac n{2}]}\cos^{2k} \varphi\cdot\sum\limits_{\substack{S\subset [1..n]\\ \#S=2k}}\prod_{j\in S}(1-a_j)\prod_{j\not\in S}a_j\le 2\sum_{k=0}^{[\frac n{2}]}\sum\limits_{\substack{S\subset [1..n]\\ \#S=2k}}\prod_{j\in S}(1-a_j)\prod_{j\not\in S}a_j.$$
The last expression is exactly $T(0)=1$, as was to be proved.
\end{Proof}

\section{The validation algorithm}
Based on the previous theorem we create an algorithm for checking the sufficient condition on C++:

\begin{lstlisting}
#include <iostream>
#include <vector>
#include <cmath>
#include <iomanip>

int main() {
    int n;
    std::cin >> n;

    std::vector<double> x(n+1, 0);
    std::vector<double> a(n+1, 0);
    std::vector<double> sigma (n+1, 1);

    for (size_t i = 1; i <= n; i ++) // compute a_i
    {
        std::cin >> x[i];
        a[i] = 1 + 2 * x[i] / ((x[i] - 1) * (x[i] - 1));
    }

    for (int k = 1; k <= n; k ++) // compute the symmetric polynomials 
    {
        double sum_j = 0;
        for (int j = 0; j <= k - 1; j ++)
        {
            double pkj = 0;
            for (int l = 1; l <= n; l ++)
            {
                pkj += std::pow(a[l], k-j);
            }
            sum_j += pow(-1, k-j-1 ) * sigma[j] * pkj;
        }
        sigma[k] = sum_j / (double) k;
    }

    std::vector<std::vector<double> > 
        c(n+1, std::vector<double>(n+1, 1));
    for(int i = 0; i <= n; ++i) // compute the binomial coefficients
    {
        for (int j = 1; j < i; ++j) 
        {
            c[i][j] = c[i - 1][j - 1] + c[i - 1][j];
        }
    }

    std::vector<double> ro(n+1, 0);
    for (int k = 0; k <= n; k++) // the symmetric means
    {
        ro[k] = sigma[k] / c[n][k];
    }

    std::vector<std::vector<double> > 
        deltaRo(n+1, std::vector<double>(n+1, 0));
    for (int i = 0; i <= n; i++) 
    {
        deltaRo[0][i] = ro[i];
    }

    for (int i = 1; i <= n; i++) // the divided differences
    {
        for (int j = i; j <= n; j++)
        {
            deltaRo[i][j] = deltaRo[i-1][j] - deltaRo[i-1][j-1];
        }
    }

    for (int i = 0; i <= n; ++i) // output
    {
        for (int j = 0; j <= n; ++j)
            std::cout << deltaRo[i][j] << " ";
        std::cout << std::endl;
    }

    bool isDone = true;
    for (int k = 0; k <= n/2; k ++)
        if (deltaRo[2 * k][n] < 0)
            isDone = false;

    if (isDone)
        std::cout << "[TRUE] The inequality holds" << std::endl;
    else
        std::cout << "[FALSE] The criteria doesn't answer" 
            << std::endl;
    return 0;
}

\end{lstlisting}
This program takes as input a natural number $n$, it is the degree of the polynomial, and then $n$  real numbers, they are the roots of the polynomial. Note that there must be at least one root equal to $-1$. The values of elementary symmetric polynomials are calculated recursively using Newton's formula: $\sigma_k=\frac 1{k}\sum_{i=0}^{k-1}\Bigl((-1)^{k-i-1}\sigma_i\sum_{j=1}^{n}a_j^{k-i}\Bigr).$ The screen displays a table of divided differences for $\rho_{k}$ and the result of checking the criterion.

\section*{Acknowledgments}
The first author is supported by the
Russian Science Foundation under grant No. 18-11-00055.

\end{document}